\documentclass{amsart}
\usepackage{graphicx,amsmath,fullpage, amssymb,color,hyperref,float}
\newtheorem{thm}{Theorem}[section]
\newtheorem{cor}[thm]{Corollary}
\newtheorem{lem}[thm]{Lemma}
\newtheorem{prop}[thm]{Proposition}
\theoremstyle{definition}
\newtheorem{defn}[thm]{Definition}

\theoremstyle{remark}

\theoremstyle{remark}
\newtheorem{rem}[thm]{Remark}
\theoremstyle{remark}

\theoremstyle{remark}
\newtheorem{conjecture}[thm]{Conjecture}
\theoremstyle{remark}
\newtheorem{convention}[thm]{Convention}
\theoremstyle{remark}

\newcommand{\To}{\longrightarrow}

\newcommand{\TD}{\mathcal{TD}}

\newcommand{\cS}{\mathcal{S}}
\newcommand{\cC}{\mathcal{C}}

\newcommand{\cN}{\mathcal{N}}
\newcommand{\cT}{\mathcal{T}}

\def\arXiv#1{{\href{http://front.math.ucdavis.edu/#1}{arXiv:#1}}}

\catcode`\@=11
\long\def\@makecaption#1#2{%
    \vskip 10pt
    \setbox\@tempboxa\hbox{
      \small\sf{\bfcaptionfont #1. }\ignorespaces #2}%
    \ifdim \wd\@tempboxa >\captionwidth {%
        \rightskip=\@captionmargin\leftskip=\@captionmargin
        \unhbox\@tempboxa\par}%
      \else
        \hbox to\hsize{\hfil\box\@tempboxa\hfil}%
    \fi}
\font\bfcaptionfont=cmssbx10 scaled \magstephalf
\newdimen\@captionmargin\@captionmargin=2\parindent
\newdimen\captionwidth\captionwidth=\hsize
\catcode`\@=12

\begin{document}
\title[]{Classification of Flat Virtual Pure Tangles}
\author{Karene Chu}
\address{Department of Mathematics at the University of Toronto\\
Fields Institute\\
Toronto Ontario \\
Canada
}
\email{karene@math.toronto.edu}

\date{\today}
\begin{abstract}

Virtual knot theory, introduced by Kauffman~\cite{Kauffman:VirtualKnotTheory}, is a generalization of classical knot theory of interest because its finite-type invariant theory is potentially a topological interpretation~\cite{Bar-Natan:UVW} of Etingof and Kazhdan's theory of quantization of Lie bi-algebras~\cite{EtingofKazhdan:BialgebrasI}.  Classical knots inject into virtual knots~\cite{Kuperberg:UInjectV}, and flat virtual knots~\cite{Manturov:FreeKnots,Manturov:FreeKnotsLinks} is the quotient of virtual knots which equates the real positive and negative crossings, and in this sense is complementary to classical knot theory within virtual knot theory.

We completely classify flat virtual tangles with no closed components (pure tangles).  This classification can be used as an invariant on virtual pure tangles and virtual braids. 

\end{abstract}

\maketitle
\newpage
\tableofcontents

\section*{Acknowledgments}
I am indebted to my PhD advisor Dror Bar-Natan for the computational evidence~\cite{BHLR:vDims} and conjecture on the main results of this paper, the idea of using the finger move in the first proof, on top of his generosity with his time, resources, care, and countless hours of inspiring discussion.  \\

 I am very grateful for the generous encouragements that other mentors have shown me, in particular Joel Kamnitzer, Louis Kauffman, and Vassily Manturov.\\

I am grateful for mathematical discussion with P. Lee, Z. Dancso, L. Leung, I. Halacheva, J. Archibald.  In particular, P. Lee pointed out~\cite{BEER:VirtualBraids} which helped me understand more about flat virtual braids, known as the ``triangular group'' in the paper.  I first learned about flat virtual knots in a lecture by Vassily Manturov in the Trieste summer school on knot theory.  \\

\section{Introduction}


We study virtual knots because they are a natural generalization of classical knots into which classical knots inject, and more interestingly, because the $R$-matrix invariants on classical knots extend naturally to virtual knots (or at least a variant of them).\\

We will define virtual knots by first recalling the definition of classical knots and generalize from it.  Classical knots can be defined combinatorially as knots diagrams modulo Reidemeister moves.  Knot diagrams are planar directed graphs with ``crossings'' as vertices.  Crossings are special tetravalent vertices whose half-edges are cyclically-ordered and directed such that opposite pairs are ``in-out'' pairs. 
The crossings have exactly the combinatorial information to be represented as follows:
\begin{figure}[H]\begin{center}
\includegraphics[height=0.5in]{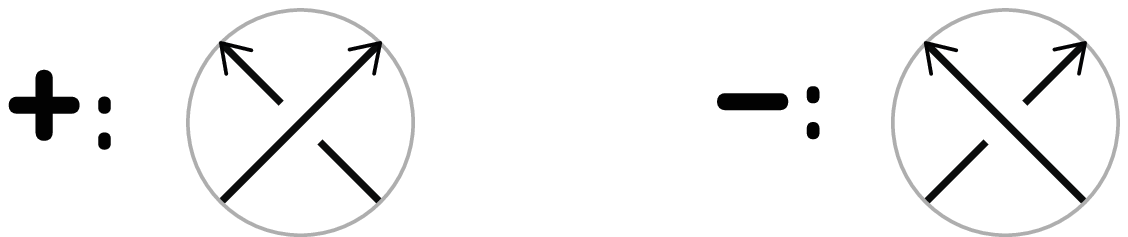}
\end{center}\end{figure}
The Reidemeister moves are local planar graph equivalence relations shown below where each skeleton strand can be oriented either way.
\begin{figure}[H]\begin{center}
\includegraphics[width=6.5in]{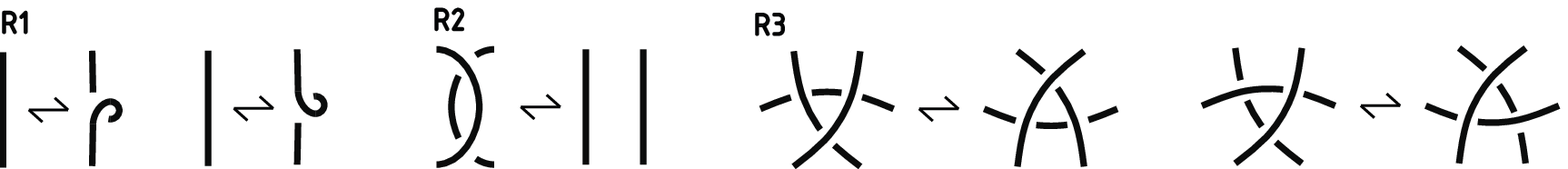}
\end{center}\end{figure}
Virtual knots have the same definition except with word ``planar'' omitted, i.e. virtual knot diagrams are (not-necessarily planar) graphs with crossings as vertices, and virtual knots are equivalence classes of virtual knot diagrams under the Reidemeister relations as local graph relations. Now, when not-necessarily planar graphs are drawn (or immersed) on the plane, transverse intersections of edges of the graph may occur.  These are not vertices of the graph, but rather artifacts of drawing a non-planar graph on the plane, and are called ``virtual crossings.''  Here is a virtual knot diagram drawn in two different ways on the plane where the real crossings are circled and the other intersections are virtual crossings:
\begin{figure}[H]\begin{center}
\label{fig:VirtualKnotDiagram}
\includegraphics[height=0.9in]{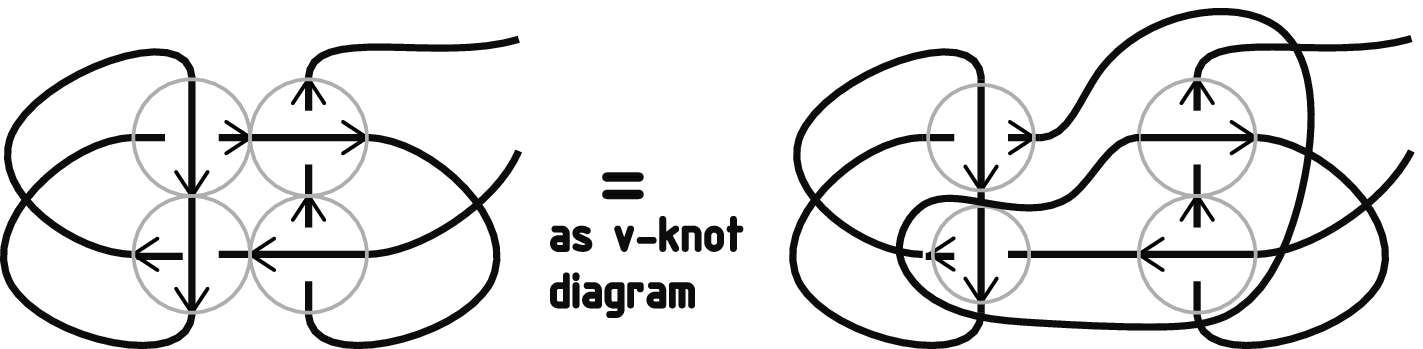}
\end{center}\end{figure}
Similarly, the strand between any two crossings in a Reidemeister relation may intersect other strands in the virtual knot diagram when drawn on the plane and have virtual crossings on them.\\

 A natural question arises: how much bigger are virtual knots than classical knots?  This leads to the consideration of the quotient of virtual knots by the crossing-flip relation which equates the (real) positive and negative crossings:
\begin{figure}[H]\begin{center}
\includegraphics[height=0.5in]{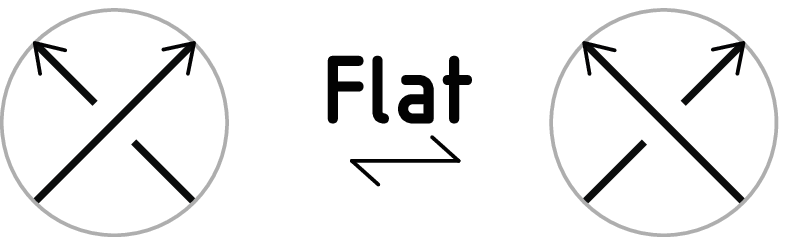}
\end{center}\end{figure}
in which all classical knots are equivalent to the unknot.  This quotient is called flat virtual knots.\\

The subject of this paper is flat virtual long knots,  where ``long'' simply refers to the skeleton of the knot being a long line, and the ``skeleton'' is the union of lines and/or circles obtained from tracing the knot diagram along the direction of the edges, across paired half-edges at crossings and forgetting the crossings.  \\

 In fact, flat virtual long knots are equivalent to descending virtual long knots, the subset of virtual long knots with only crossings whose over strand is earlier w.r.t. to the orientation of the skeleton than the under strand.  This is because while virtual long knots project onto flat virtual long knots, there is a well-defined section map from flat virtual long knot back into virtual long knot, namely by sending any flat real crossing to a descending crossing,  and the image of such a section map is exactly the descending virtual long knots.  Notice descendingness of a crossing is not defined for round virtual knots, and the map that sends any flat crossing to a positive crossing is well-defined.\\

Our main result is the classification of both the ``framed'' and ``unframed'' versions of descending virtual long knots, where ``framed'' means the Reidemeister $1$ relation is not imposed and ``unframed'' means otherwise.  We give a canonical representative for each equivalent class of descending virtual long knot diagrams under the Reidemeister moves.

\begin{figure}[H]\begin{center}
\includegraphics[height=1.3in]{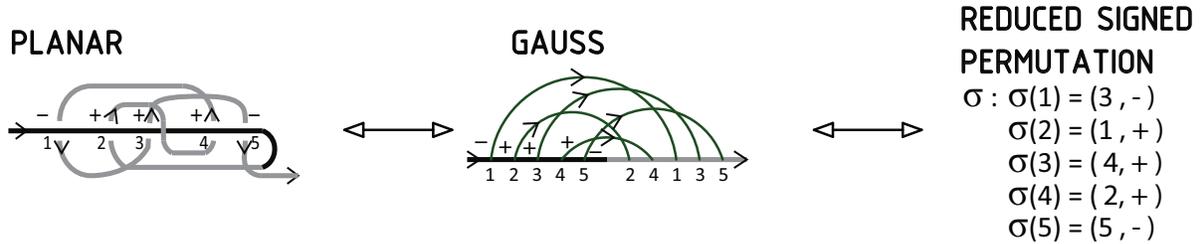}
\caption{Example of the canonical form of a descending virtual long knot.  There is a point on the skeleton before which it is the over strand in all crossings it participate in and after which it is under.}
\label{fig:CanonicalKnotExample}
\end{center}\end{figure}

\begin{figure}[H]\begin{center}
\includegraphics[width=6.2in]{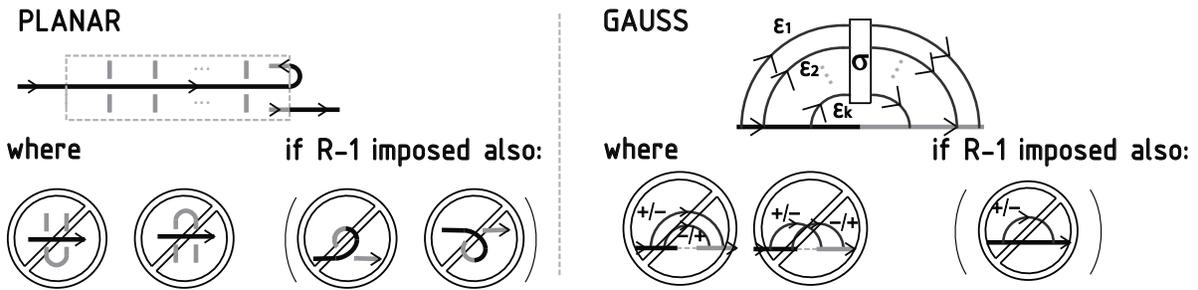}
\caption{General form of the \textbf{canonical diagrams} of descending virtual long knots, characterized by the existence of a point on the skeleton before which it is the over strand in all crossings it participates in and after which it is under, and the exclusion of the bigons and ``R1-kinks'' as well for the unframed version.  In the Gauss diagrams on the right, the $\epsilon$'s are signs of crossings, and the box with $\sigma$ denotes a permutation of the arrows so that the incoming arrows are permuted by $\sigma$ within the box and emerge on the other side permuted. }
\label{fig:CanonicalKnots}
\end{center}\end{figure}

\begin{figure}[H]\begin{center}
\label{fig:KnotTable}
\includegraphics[width=6in]{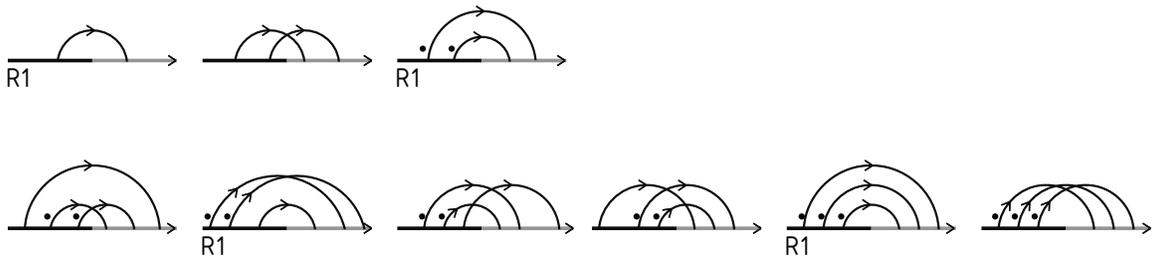}
\caption{List of canonical diagrams in Gauss diagram form of framed flat virtual long knots up to three crossings.   For the unframed version, exclude the diagrams with $R1$ below them.  Chords in the same diagram with dots on their left are required to have the same signs; all others can be either $+$ or $-$. Thus, the first diagram in the second row represents $2\times2$ different canonical diagrams with different sign arrangements.}
\end{center}\end{figure}

\begin{thm}[Classification of Long Descending Virtual Knots, conjectured by Bar-Natan]
\label{thm:ClassifyFlatKnots}
 Framed descending virtual long knots $\mathcal{K}^{\text{\it{f}}}$ are in bijection with the set of canonical diagrams $\cC^{1}$.   A \textbf{canonical diagram} is a descending virtual long knot diagram whose skeleton strand has a point before which it is the over strand in any crossing it participates in, and after which as the under strand, and does not contain bigons bounded by opposite signed crossings.  An example is given in figure~\ref{fig:CanonicalKnotExample} and the general form is shown in figure~\ref{fig:CanonicalKnots}.   \\Furthermore, $\cC^{1}$ is in bijection with the set of all ``signed reduced permutations'', where a \textbf{signed permutation} is a set map $\tilde{\rho}: \{1,\ldots, n\}\longrightarrow \{1,\ldots, n\} \times \{+, - \}$ which projects to the first components as a permutation, and a \textbf{reduced signed permutation} satisfies the extra condition that the image of pairs of consecutive numbers are not pairs of consecutive numbers with opposite signs, i.e. not $((j,\mp), (j+1,\pm))$, or $((j+1,\pm),(j,\mp))$ for all $j<n$.\\
Long unframed flat virtual knots $\mathcal{K}^{f}$ are in bijection with the subset of $\cC^{1}$ with no ``R1 kinks,'' also shown in ~\ref{fig:CanonicalKnots}.  See figure~\ref{fig:KnotTable} for a list of canonical diagrams up to three crossings.
\end{thm}

\begin{figure}[h!]\begin{center}
  \includegraphics[height=1.5in]{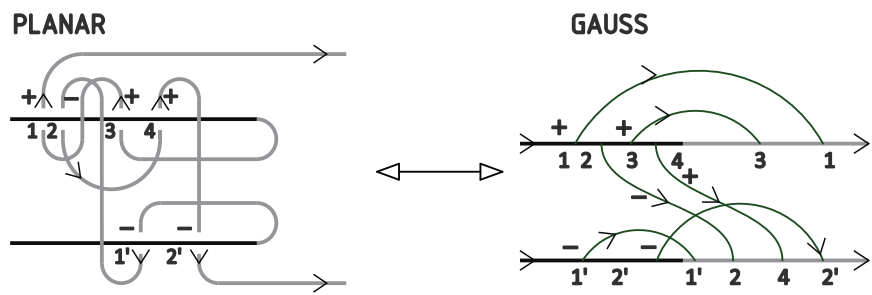}
 \caption{The canonical diagram of a framed descending virtual pure tangle on two strands.  On each skeleton strand, there is a point before which it is over in all crossings and after which it is under in all crossings, and the diagram contains no bigons bounded by opposite signed crossings.  Since it also does not contain any R1-kinks, it is a canonical diagram also of an unframed descending virtual pure tangle.}
 \label{fig:CanonicalTangleExample}
 \end{center}\end{figure}

\begin{figure}[H]\begin{center}
\includegraphics[width=6in]{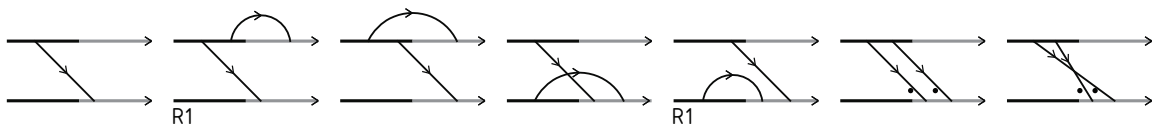}
\caption{List of canonical diagrams in Gauss diagram form up to two crossings of framed descending virtual pure tangle on two strands with at least one crossing between the two strands.  For the unframed version, exclude the diagrams with $R1$ below them.  Notice the top and bottom strands are distinguishable since the strands are ordered.  All dotted chords in the same diagram are required to have the same signs, and all other are signed in all ways possible.}
 \label{fig:TangleTable}
 \end{center}\end{figure}

The above result can be generalized easily to the multi-strand case.  We call flat virtual tangles whose skeleton is an ordered union of strands (in particular no closed loops) flat virtual pure tangles.  Similar to the long knot case, these are again equivalent to descending virtual pure tangles.
\begin{thm}[Classification of Descending Virtual Pure Tangles]
\label{thm:ClassifyFlatTangles}
Framed descending virtual pure tangles of $n$ strands $\mathcal{T}^{\text{\it{f}}}$ are in bijection with the set of canonical diagrams $\cC^{n}$, which are characterized by the same two conditions as in the one strand case in theorem~\ref{thm:ClassifyFlatKnots} but applied to all $n$ strands.   Unframed descending virtual pure tangles are in bijection with the subset of $\cC^{n}$ with no ``R-I kinks''.  See figure~\ref{fig:CanonicalTangleExample} for an example, and figure~\ref{fig:TangleTable} for a partial list of canonical diagrams of descending virtual pure tangles on two strands up to two crossings.
\end{thm}%
%

The proofs of both theorems are similar and amount to showing that a well-defined sorting map exists.  For an example of the sorting, see figure~\ref{fig:SortExample}.\\

This classification can be use as an invariant on virtual long knots and pure tangles, as well as virtual braids.  If flat virtual braids inject into the pure tangles, then we have also obtained their classification.\\

\section{Classification of Pure Descending Virtual Tangles}
\label{sec:Global}

Having established that flat virtual pure tangles are equivalent to descending virtual pure tangles in the introduction,
we present in this section the classification of descending virtual pure tangles and its proof.   We will prove the classification for the framed version (see page~\pageref{thm:ClassifyFlatKnots}) first, and then modify it slight for the unframed version to include Reidemeister 1. 

\subsection{Generic Diagrams of Pure Descending Virtual Tangles}

In this subsection, we describe the generic form of pure descending virtual tangle diagrams.  First, a few definitions to describe the diagrams:
\begin{defn}
An interval of the skeleton of a pure descending virtual tangle is called an \emph{over} (resp. \emph{under}) interval if all of its subintervals that take part in crossings are the over strands in the crossings.  A \emph{maximal over} (resp. \emph{maximal under}) interval is an over (resp. under) interval preceded and followed immediately by an under (resp. over) interval or by the beginning or end of the strand.  An \emph{illegal interval} is an interval consisting of first a maximal under interval and then a maximal over interval.  These are illustrated below in Figure~\ref{fig:IllegalInterval}.
\end{defn}

 \begin{convention}
All real crossings are descending all virtual crossings are neither over or under and not circled.  The over interval of a crossing is drawn in black and the under interval drawn in grey.  An interval not explicitly oriented means it can be oriented either ways.  In any Gauss diagram, an unsigned Gauss arrow means it can have either sign.  A ``thick band'' represents multiple strands or arrows, as in figure~\ref{fig:IllegalInterval}.
 \end{convention}

\begin{figure}[H]\begin{center}
\includegraphics[height=0.7in]{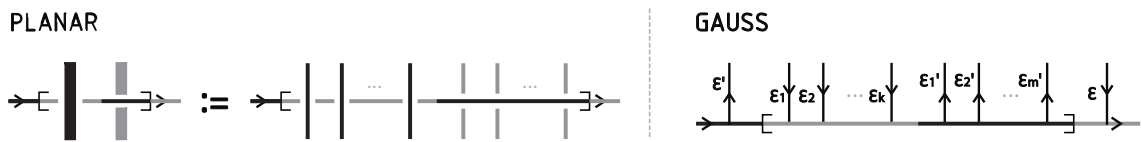}
\caption{An \textbf{illegal interval}, within the square brackets, is  the skeleton interval which is first a maximal under interval (grey) and then a maximal over interval (black).  Any subintervals of the maximal under (resp. over) is an under (resp. over) interval.  The interval preceding (resp. following) this illegal interval is either an over (\emph{resp.} under) interval or the beginning (\emph{resp.} end) of the skeleton strand.  In the Gauss diagram on the left the illegal interval is one which is in between an over and an under interval.   The half arrows have their other ends on other parts of the skeleton. }
\label{fig:IllegalInterval}
\end{center}\end{figure}

\begin{figure}[H]\begin{center}
\includegraphics[height=1.5in]{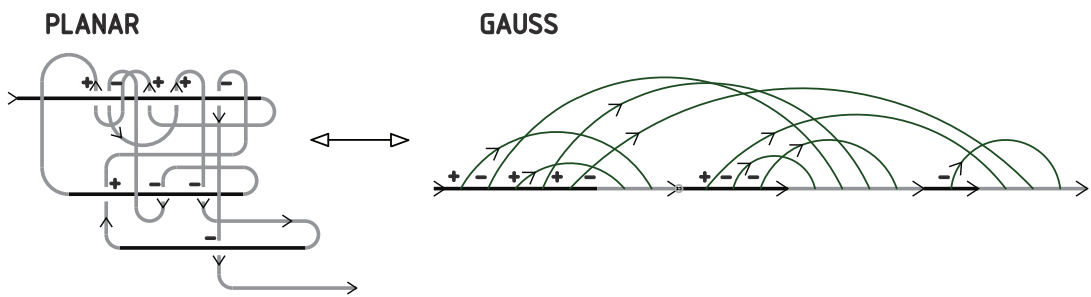}
\caption{A generic diagram for a descending virtual long knot.  The skeleton strand can be partitioned into maximal over and under intervals.  It starts with a maximal over one, then alternates between maximal over and under, and ends in a maximal under interval.  The maximal over interval are in black, and the under in grey in both the planar and Gauss diagrams.}
\label{fig:GenericKnot}
\end{center}\end{figure}

Generically, a pure descending virtual long knot diagram has multiple maximal over and under intervals.  Due to descendingness, it always (as long as there is at least one crossing) starts with a maximal over and ends with a maximal under interval, while in between it alternates between over and under while having each maximal under interval only under the maximal over intervals before it.  Thus, a generic diagram has illegal intervals on its skeleton. See figure~\ref{fig:GenericKnot} for an example.  A generic descending virtual pure tangle diagram is simply a descending virtual long knot diagram with a finite number of cuts on its skeleton.

%

\begin{rem}
\label{rem:TDParameters}
 There are two parameters on the set of descending virtual pure tangle diagrams:  the number of illegal intervals, $\cN(D)$, and the number of crossings, $\chi(D)$ of a diagram $D$.  Both are non-negative for all diagrams.  Furthermore, the number of crossings is bounded below by $\chi(D)\geq N(D)+1$, since in the Gauss diagram language, each of $\cN(D)$ illegal intervals in a diagram $D$ must have at least one arrow-head and one arrow-tail, summing to $2N$ half arrows within the illegal interval, and the beginning of the first strand and the end of the last strand must have one arrow-tail and one arrow-head respectively. And this bound is attained by the following diagram:
\begin{figure}[H]\begin{center}
\includegraphics[width=3in]{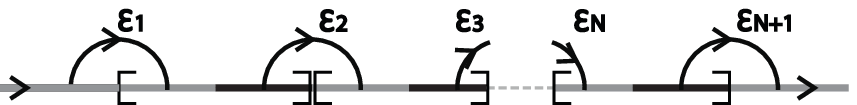}
\end{center}\end{figure}
\end{rem}
\subsection{The Sorting Map}
\label{sec:SortingMap}
We start presenting the proofs of theorems~\ref{thm:ClassifyFlatKnots} and~\ref{thm:ClassifyFlatTangles}.  Refer to page~\pageref{thm:ClassifyFlatKnots} for the definition of \emph{framed} and \emph{unframed}, statement of the theorems and in it the definitions of \emph{canonical diagrams}, and \emph{reduced signed permutations}.

We first show the bijection (in theorem~\ref{thm:ClassifyFlatKnots}) between the canonical diagrams $\cC^{1}$ for framed long descending virtual knots and reduced signed permutations, and then describe a sorting map $\cS$ that chooses a canonical representative diagram for each class of equivalent framed descending virtual pure tangle diagrams. \\

\begin{prop}
 The set of canonical diagrams $\cC^{1}$ (figure~\ref{fig:CanonicalKnots}) of framed descending virtual long knot is in bijection with the set of reduced signed permutations.
\begin{proof}
Consider a canonical diagram $C$ with $n$ arrows in the Gauss Diagram language.  Label the arrow-tails by $1,2,\ldots,n$ in increasing order from the start of the knot, and label the arrow-heads similarly beginning with the first arrow head.  Then construct a reduced signed permutation $\tilde{\rho}$ from the diagram by $\tilde{\rho}(i)=(j,\epsilon)$ where $j$ and $\epsilon$ are respectively the arrow-head label and the sign of the arrow with tail labeled $i$.  There being no bigons in the canonical diagram translates to the restriction that the image under $\tilde{\rho}$ of pairs of consecutive numbers are not any of $((j,\mp), (j+1,\pm))$, and $((j+1,\pm),(j,\mp))$ for any $j<n$.  The inverse of this map is obvious.
\end{proof}
\end{prop}

First, we introduce the finger move, \emph{F-move}:
\begin{figure}[H]\begin{center}
\includegraphics[width=6in]{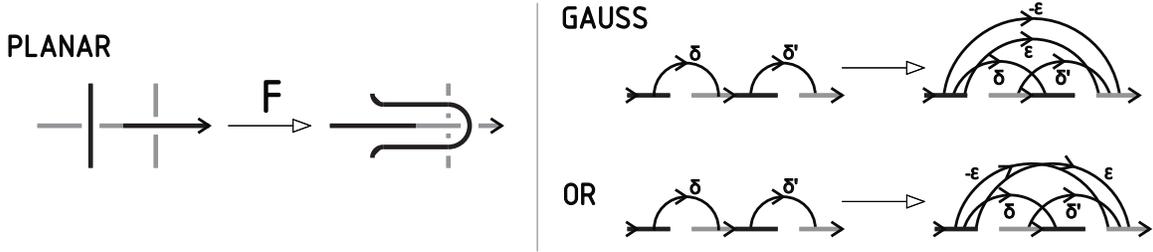}
\caption{\textbf{Finger move}.  In the Gauss diagram language, there are two resulting diagrams depending on the relative orientations of the two vertical strands in the planar diagram.  $\delta$'s and $\epsilon$'s are signs. }
\end{center}\end{figure}

\begin{prop} The set of all $F$-moves and $R2$-moves is equivalent to the set of all $R3$-moves and $R2$-moves.
\label{prop:FMovesR2}
\begin{proof}
$R3$-moves are generated by $R2$-moves and $F$-moves, as shown in the following figure. Similarly, $F$-moves are generated by $R2$- and $R3$-moves:
\begin{figure}[H]\begin{center}
\includegraphics[height=1.3in]{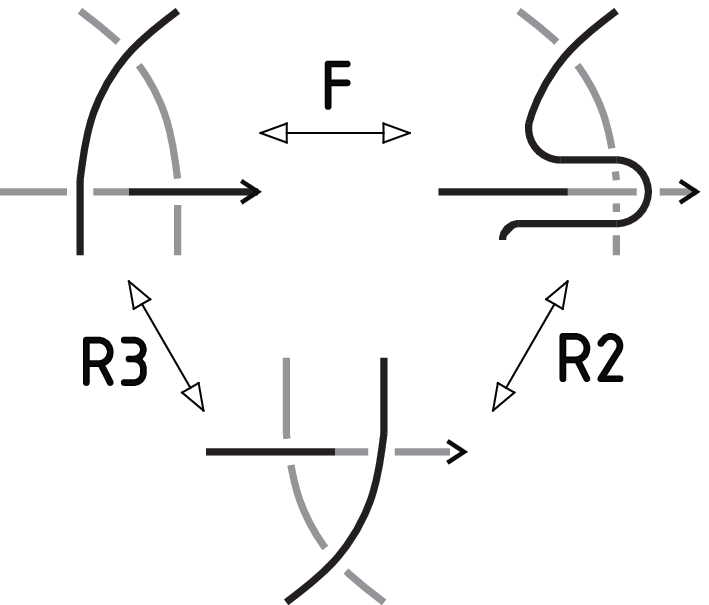}
\end{center}\end{figure}
Half of the R3-moves are represented by the above diagram, the other half are represented by the up-down-mirror image of the diagram.
\end{proof}
\end{prop}

\begin{cor}
To show that a map on the diagrams of descending virtual pure tangles descends to the equivalence classes of framed descending virtual pure tangles, it suffices to show that the map is well defined under the finger moves and R-2 moves.
\end{cor}

From now on we omit the Gauss diagrams since they have become too complicated to draw but they can be constructed easily. 

 The two following local sorting moves will be used in the sorting of a generic diagram into its canonical form.
\begin{defn}
The sorting group-finger-move, \emph{GF-sort}, and the sorting R2-move, \emph{R2-sort}, are the following single-direction moves that take place inside the squared region, called the \emph{sorting site}:
\begin{figure}[h!]\begin{center}
\includegraphics[width=5.5in]{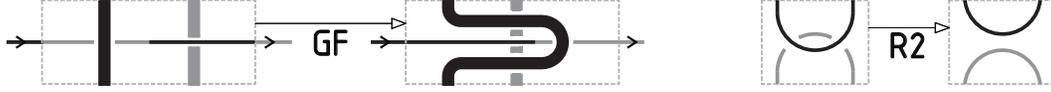}
\caption{(L) GF-sort;  (R) R2-sort.  An over (\emph{resp.} under) thick band denotes multiple over (\emph{resp.} under) strands, as shown in figure~\ref{fig:IllegalInterval} before.  Notice these sorting moves go only in one direction.}
\end{center}\end{figure}
\end{defn}

\begin{rem}
\label{rem:R2GFParameters}
 \begin{enumerate}
          \item GF-sort is generated by single sorting F-moves and so is generated by R2 and R3-moves.
          \item GF-sort switches the order of the maximal over interval and maximal under intervals within the illegal interval, thus decreasing the number of total illegal intervals by 1, but lengthens the preceding maximal over and the following maximal under intervals.
          \item GF-sort increases the number of total crossings by $2n>0$ of the diagram.
          \item R2-sort decreases the number of total crossings by $2$, and either does not change or decreases the number of total illegal intervals by at most $2$.
        \end{enumerate}
\end{rem}

Some more terminology for the definition of the sorting map.
\begin{defn}
\noindent
\begin{enumerate}
\item
A sorting move is \emph{available} in a diagram $D$ if a subdiagram of $D$ is equal to the L.H.S. of the sorting move.  This subdiagram is called the \emph{sorting site} in $D$ for the sorting move;\\
\item
Two sorting moves $s$, $t$ \emph{overlap} if in the intersection of their sorting sites, there is at least a crossing.\\
\item
A \emph{sort sequence} $S$ on a diagram $D$ is a finite sequence of sorting moves $s_k \circ \ldots\circ  s_2 \circ s_1$ such that for each $i$, $s_i$ is an available move on the diagram $s_{i-1} \circ  \ldots\circ  s_2 \circ s_1(D)$.
\item
A \emph{terminating sort sequence} on $D$ is a sort sequence $T$ such that $T(D)$ has no available sorting moves.
\end{enumerate}
\end{defn}

We can now characterize the set of canonical diagrams $\cC$ to be all framed descending virtual pure tangle diagrams with no illegal intervals and no $R2$-sorting sites.

\begin{defn}
\label{defn:SortingMap}
Define the \emph{sorting map} on the set of all pure descending virtual tangle diagrams $\mathcal{TD}^{\text{\it{f}}}$ to be
\begin{eqnarray}
\nonumber \mathcal{S}: \mathcal{TD}^{\text{\it{f}}} &\To& \mathcal{TD}^{\text{\it{f}}} \\
\nonumber D &\longmapsto& s_k \circ\ldots \circ s_2\circ s_1(D)
\end{eqnarray}
where $s_k\circ \ldots\circ s_2 \circ s_1$ is any terminating sort sequence on $D$.  We show below that $\mathcal{S}$ is well-defined.
\end{defn}

See section~\ref{sec:Examples} for examples of sorting.

\begin{prop}
\noindent
\begin{enumerate}
\item \label{item:SGen} $\cS$ is generated by Reidemeister-moves;
\item $\cS$ is defined, i.e. the algorithm terminates
\item For any pure tangle diagram $D$, $\cS(D)\in \cC \in\TD^{\text{\it{f}}}$
\end{enumerate}
\begin{proof}
\begin{enumerate}
\item Both GF- and R2 sorts are a finite sequence of Reidemeister-moves;
\item Only finite number of GF-sorts can be performed since a GF-sort decreases the parameter $N_{D}$ (the number of illegal intervals) by 1 and R2-sorts do not increase $N_D$.  Since the number of GF-sorts are finite, at the point in any sorting algorithm when all GF-sorts are performed, only finite R2-sorts can be performed since it decreases the parameter $\chi_{D}$ by 2;
\item The result of any terminating sort sequence has no illegal intervals and no R2-sorting sites.
\end{enumerate}
\end{proof}
\end{prop}

\begin{lem}
\label{lem:SWellDefined}
$\cS:\TD^{\text{\it{f}}} \To\TD^{\text{\it{f}}}$ descends to a bijection $\cS:\cT^{\text{\it{f}}} \To\cT^{\text{\it{f}}}$ between the set pure descending virtual tangles $\mathcal{T}^{\text{\it{f}}}$ and the set of canonical diagrams $\cC$ (defined in theorems~\ref{thm:ClassifyFlatKnots} and~\ref{thm:ClassifyFlatTangles} on page~\pageref{thm:ClassifyFlatKnots}.)
\begin{proof}
We need to show that $\cS$ is well-defined under choices of terminating sorting sequences, and well-defined under Reidemeister-moves, and is bijective into the set of canonical diagrams $\cC$.  Well-definedness of $\cS$ follows from lemmas~\ref{lem:SWellDefinedSorts} and~\ref{lem:SWellDefinedRMoves} in the next section.  It remains to show bijectivity, but surjectivity follows from the fact that a canonical diagram does represent a pure descending virtual tangle and injectivity follows from the fact that $\cS$ applied to any canonical diagram results in the same canonical diagram.
\end{proof}
\end{lem}


\subsection{Sorting map is well-defined}
\label{sec:SortingMapWellDefined}
This section is the main part of the proof of lemma~\ref{lem:SWellDefined}, divided into lemmas~\ref{lem:SWellDefinedSorts} and~\ref{lem:SWellDefinedRMoves}

\begin{lem}  $\mathcal{S}$ is well-defined under choices of different terminating sort sequences.  \label{lem:SWellDefinedSorts}
\begin{proof}
We proceed by a two-dimensional induction on $(\cN(D),\chi(D))$, the number of illegal intervals and the number of crossings of a diagram  $D\in \TD^{\text{\it{f}}}$.  The induction steps will involve the diamond lemma.
We will first show that $\cS(D)$ is well-defined for all the diagrams $D$ in the ``column'' $\cN=0$ using an induction on the variable $\chi$, and then assuming the induction hypopaper for all ``columns'' $\cN(D)\leq n$, show the statement for the ``column'' $\cN=n+1$ by another an induction on $\chi$.  In all induction steps below, we will use one of the two following general arguments.  
We will call a region in the inductive domain where the statement is already true, either by hypopaper or by proof, a \emph{truth region}.  First, for the case when a diagram $D$ has only $1$ available sorting move, $s$, the sorting move $s$ on $D$ will result in a diagram in a truth region, i.e. any terminating sorting sequence on $s(D)$ gives the same resulting diagram.  This then implies that any terminating sorting sequence on $D$ itself results in the same diagram.  Second, for the case when a diagram $D$ has two or more available sorting moves, it suffices to show that for any pair of available sorting moves $s$ and $t$ on $D$, any terminating sort sequence starting with $s$ will give the same resulting diagram as any terminating sort sequence starting with $t$.  As in the previous case, both $s(D)$ and $t(D)$ will be in a truth region, ie. all terminating sort sequences $S$ on $s(D)$ will result in the same diagram, and the same for $t(D)$. In particular, if we can choose sorting sequences $S$ on $s(D)$ and $T$ on $t(D)$ such that $S(s(D))=T(t(D))$, the claim follows. There are two cases: if $s$ and $t$ do not overlap, they commute and the trivial relation between relations, also known as a \textbf{syzygy}, $st=ts$, can be used; otherwise, syzygies $S\circ s= T\circ t$ will be needed for the argument.\\
%
%


Thus, for all induction steps below, we only need to verify that for the given diagram $D$, any available sorting move on it does result in a diagram in the true region, and that for any pair of available overlapping sorting moves $s$,$t$ on $D$, there are specific syzygies $S(s)=T(t(D))$ to substitute into the above argument. \\

We proceed to check these for all steps in our two dimensional induction.  Also recall a sorting move is either an $R2$- or a $GF$-sort.
First, for the \textbf{Base ``column, $\cN=0$ }, any diagram with zero illegal intervals has no available $GF$-sorts.

\begin{description}
\item[\textbf{Base case, $(\cN,\chi)=(0,0)$ or $(0,1)$}]  With less than two crossings, a diagram has no available $R2$-sort either, so $\cS(D)=D$ is well-defined.
\item[\textbf{Induction, ``$\chi \leq c-1$'' $\Rightarrow$ ``$\chi=c$''}]  Assume $\cS$ is well-defined on all diagrams with $\chi\leq c-1$ where $c\geq 2$. The only possible available sorting moves on a diagram $D$ with $(\cN, \chi)=(0,c)$ are R2-sorts, and by remark~\ref{rem:R2GFParameters}, any R2-sort on $D$ will result in a diagram in the truth region ``$\chi \leq c-2$.''  Also, up to orientation of the strands, there are only two possible ways R2-sorts can overlap, and they are mirror images of one another.  The syzygy for one of these overlapping R2-sorts are shown in figure~\ref{fig:R2R2Syzygy}; the other one is analogous.
\end{description}
\pagebreak
 Secondly, the \textbf{Induction step ``columns $\cN\leq N-1$'' $\Rightarrow $ ``column $\cN=N$''}.
 \begin{description}
 \item[\textbf{Base case, ``columns $\cN \leq N-1$'' $\Rightarrow$ ``$(\cN,\chi)=(N,\chi_{\text{min}}(N))$ ''}]
  A diagram $D$ with the minimum number of crossing $\chi_{\text{min}}$ to make $N\geq1$ illegal intervals has only one crossing in each maximal over or under interval, and so has no available $R2$-sorts.  But $GF$-sorts can be available and by remark~\ref{rem:R2GFParameters} result in diagrams in a truth region.  Now, there is only one way two $GF$-sorts can overlap and there is a syzygy between them, shown in figure~\ref{fig:GFGFSyzygy}.
  \begin{figure}[h!]\begin{center}
\includegraphics[width=2.5in]{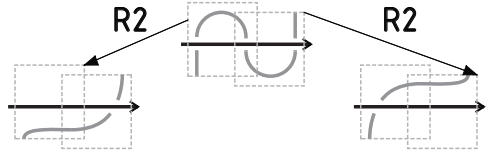}
\caption{Syzygy for overlapping $R2$-sorts: performing either available $R2$-sorts leads to the same diagram with fewer crossings. Each sorting move happens inside the corresponding sorting sites, boxed by light dotted lines, in different diagrams. }
\label{fig:R2R2Syzygy}
\end{center}\end{figure}
\begin{figure}[h!]\begin{center}
\includegraphics[width=6.5in]{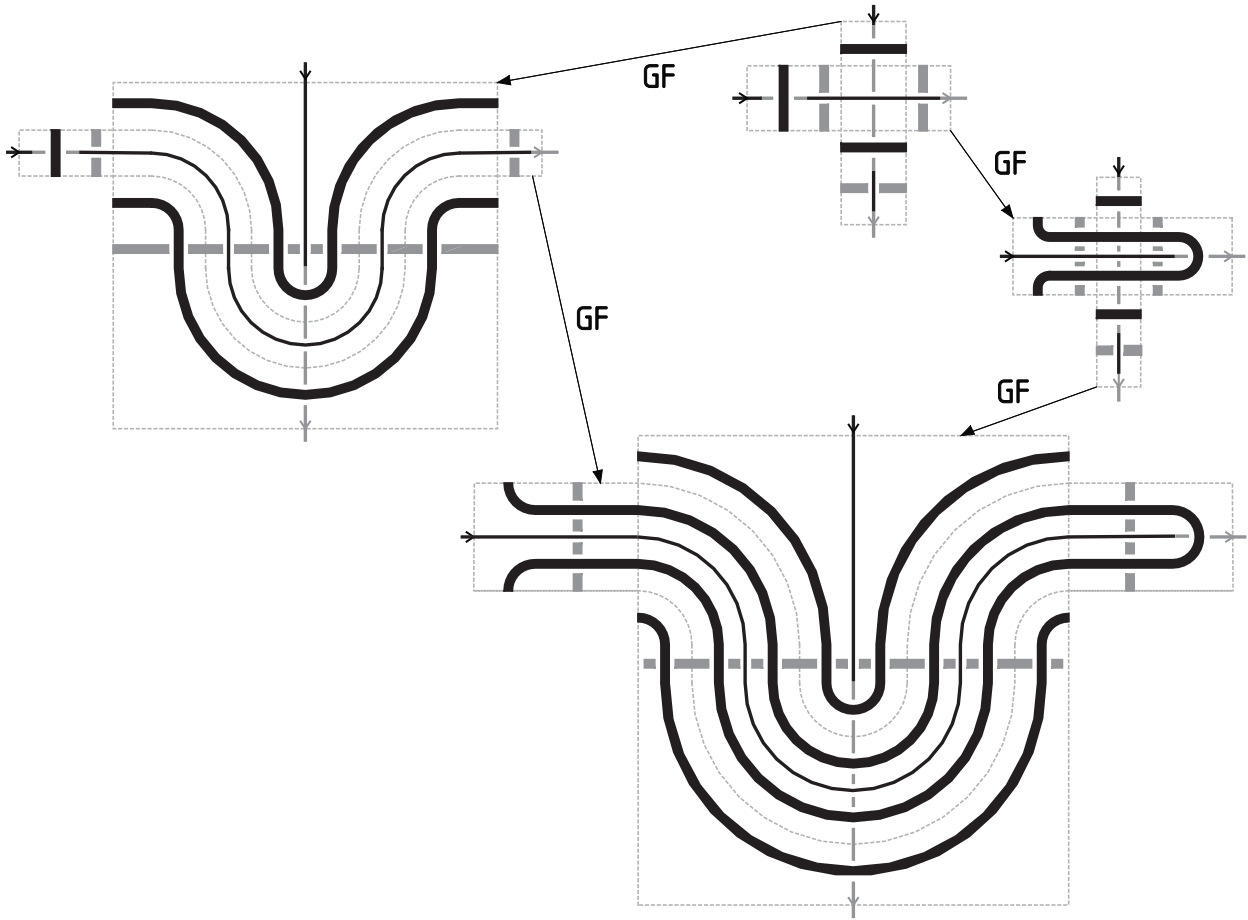}
\caption{Syzygy for overlapping $GF$-sorts.  Thick bands are multiple strands, as in figure~\ref{fig:IllegalInterval}.  There are two available GF-sorts to perform on the top diagram, with their sorting sites in the vertical and horizontal boxes respectively.  The path leading from the top diagram first to the left has the GF-sort in the vertical box performed first, followed by the only remaining GF-sort available, inside the horizontal box with a U-shape.  The path leading from the top first to the right has the GF-sort in the horizontal box performed first, followed by its only remaining available GF-sort, inside the vertical box.  Each GF-sort lowers the number of illegal interval and both paths lead to the same diagram at the bottom.}
\label{fig:GFGFSyzygy}
\end{center}\end{figure}
\item[\textbf{Induction,  columns ``$\cN \leq N-1$'' and ``$\cN=N, \chi\leq  c-1$'' and   $\Rightarrow$  ``$(\cN,\chi)=(N,c)$''}]
  Assume $\cS$ is well-defined on all diagrams with less than $N$ illegal intervals, where $N\geq 1$, and all diagrams with $N$ illegal intervals and less than $c$ crossings where $c > \chi_{min}(N)$.  Now, on a diagram $D$ with $(\cN,\chi)=(N, c)$ both $GF$- and $R2$- sorts can be available and by remark~\ref{rem:R2GFParameters}, both will result in a diagram in a truth region.  
We also need syzygies for all ways of overlap of all sorting moves, $R2$-$R2$, $GF$-$GF$, and $GF$-$R2$.  The first two cases $R2$-$R2$ and $GF$-$GF$ are the same as in previous steps, with syzygies shown in figures~\ref{fig:R2R2Syzygy} and ~\ref{fig:GFGFSyzygy}.  For the third case, a $GF$-sort and an $R2$-sort can overlap in essentially two ways up to orientation of strands, depending on whether the crossings in the $R2$-sorting site belong to the maximal over or under interval in the $GF$-sorting site.  Also, within each of these overlap types, the $R2$-sort site can still vary.  The syzygy for the first way is shown in figure~\ref{fig:GFR2Syzygy}; the one for the second is analogous.
\begin{figure}[h!]\begin{center}
\includegraphics[width=5.5in]{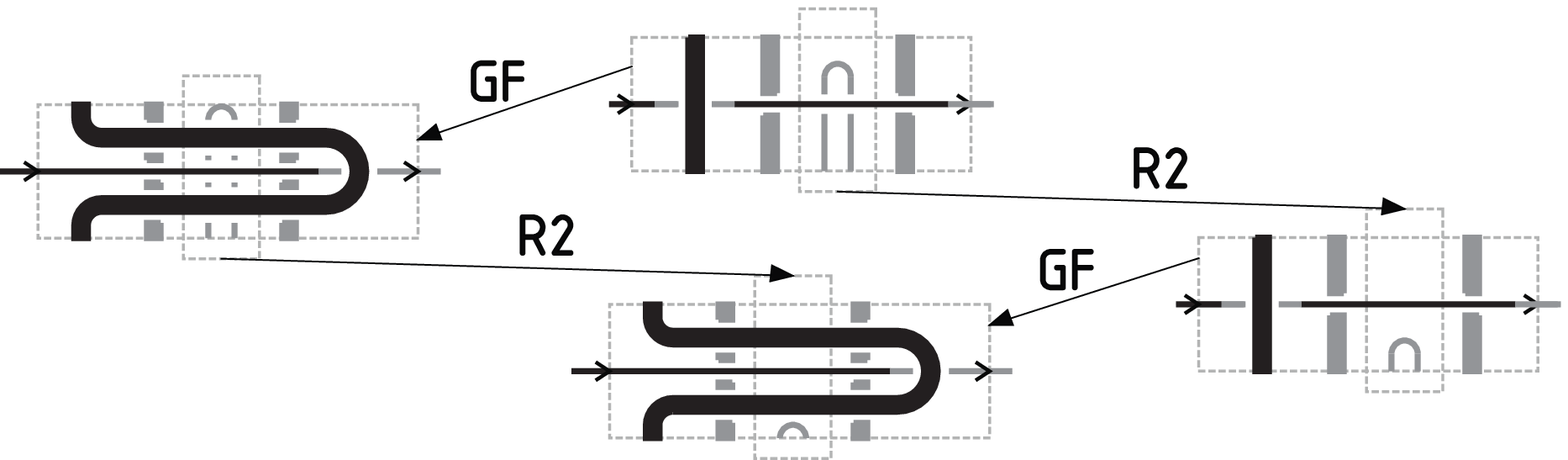}
\caption{Syzygy for overlapping $GF$ and $R2$-sorts.  Thick bands are multiple strands, as in figure~\ref{fig:IllegalInterval}.  There are two available sorts to perform on the top diagram, a $GF$-sort and an $R2$-sort. The two paths leading from the top to the bottom diagram corresponds to different choices of which of $GF$- and $R2$-sorts to perform first, and result in the same diagram with fewer illegal intervals.}
\label{fig:GFR2Syzygy}
\end{center}\end{figure}
\end{description}
\end{proof}
\end{lem}

\begin{lem}
\label{lem:SWellDefinedRMoves}
 The sorting map $\cS$ is well-defined under Reidemeister-II and III moves.
\begin{proof}
This follows directly from proposition ~\ref{prop:FMovesR2} and the next lemma ~\ref{lem:SWellDefinedFMoves}.
\end{proof}
\end{lem}

\begin{figure}[H]\begin{center}
\includegraphics[width=6.5in]{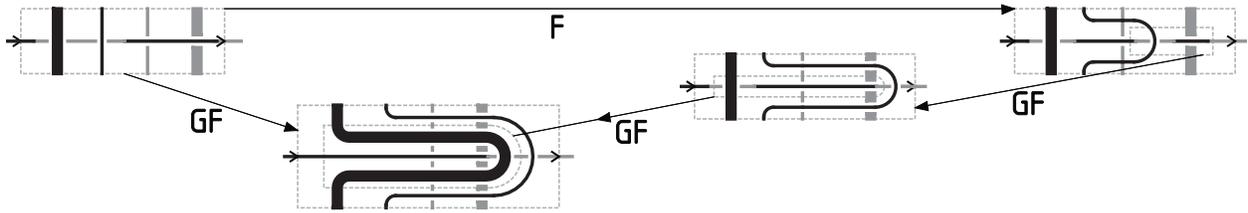}
\caption{Syzygy for overlapping $GF$-sort and finger move.  Two diagrams (on the left and right most) differing by a single $F$-move can be sorted by available $GF$-sorts to the same diagram (at the bottom).}
\label{fig:GFFSyzygy}
\end{center}\end{figure}

\begin{lem}
\label{lem:SWellDefinedFMoves}
The sorting map $\cS$ is well-defined under finger-moves.
\begin{proof}
Since $\cS$ is well-defined under choices of different terminating sort sequences on all diagrams $D\in\TD^{\text{vf}}$, if we can choose a sorting sequences on both sides $D_{l}$ and $D_{r}$ of the finger-move such that they result in the same diagram, then $\cS(D_{l})=\cS(D_{r})$.  Thus, the syzygy in figure~\ref{fig:GFFSyzygy} suffices.
\end{proof}
\end{lem}
This completes our proof of the classification of the framed version of pure descending virtual tangles, the first statements in theorems~\ref{thm:ClassifyFlatKnots},~\ref{thm:ClassifyFlatTangles}.

\pagebreak
\subsection{Classification of the Unframed Version: Adding Reidemeister I}

To prove the second statements of theorems~\ref{thm:ClassifyFlatKnots},~\ref{thm:ClassifyFlatTangles} which classify the unframed version of pure descending virtual tangles, which recall is the quotient of the framed version by the Reidemeister-I relation, we only need to slightly modify the proof of the framed version in the last sections~\ref{sec:SortingMap},\ref{sec:SortingMapWellDefined}.  First, we add an extra sorting move, the $R1$-sort as shown below in figure~\ref{fig:R1Sort},  to the definition~(\ref{defn:SortingMap}) of the sorting map.
\begin{figure}[H]\begin{center}
\includegraphics[width=6.2in]{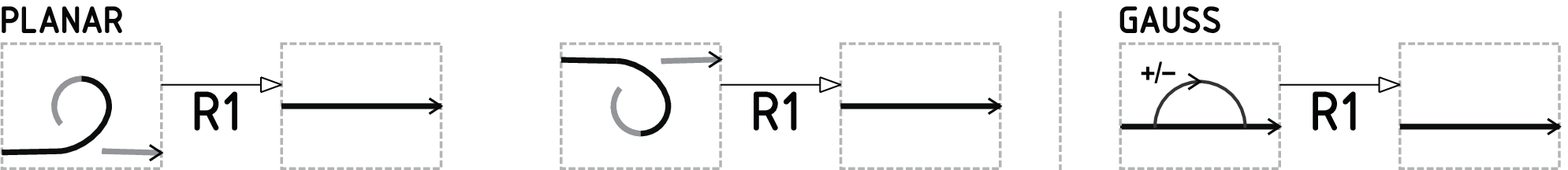}
\caption{$R1$-sort.}
\label{fig:R1Sort}
\end{center}\end{figure}

 Then, we show that the modified sorting map is still well-defined by adding R1-sort to the two-dimensional induction argument in~\ref{lem:SWellDefined}: we note that performing any R1-sort will either decrease the number of illegal intervals $\cN$ by 1 or not change it, and will always decrease the number of crossings $\chi$ by 1, thus resulting in a diagram in the already true region in the induction domain; and use the following two overlapping syzygies to conclude that the choice to perform an $R1$-sort, an $R2$-sort, or a $GF$-sort at each stage of the sorting does not affect the result.

\begin{figure}[H]\begin{center}
\includegraphics[width=6.5in]{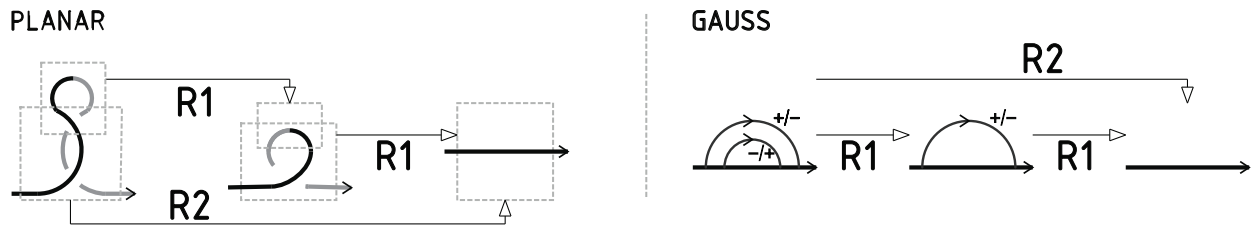}
\caption{Syzygy between $R1$- and $R2$- sorting moves.}
\label{fig:R1R2Syzygy}
\end{center}\end{figure}

\begin{figure}[H]\begin{center}
\includegraphics[width=6.5in]{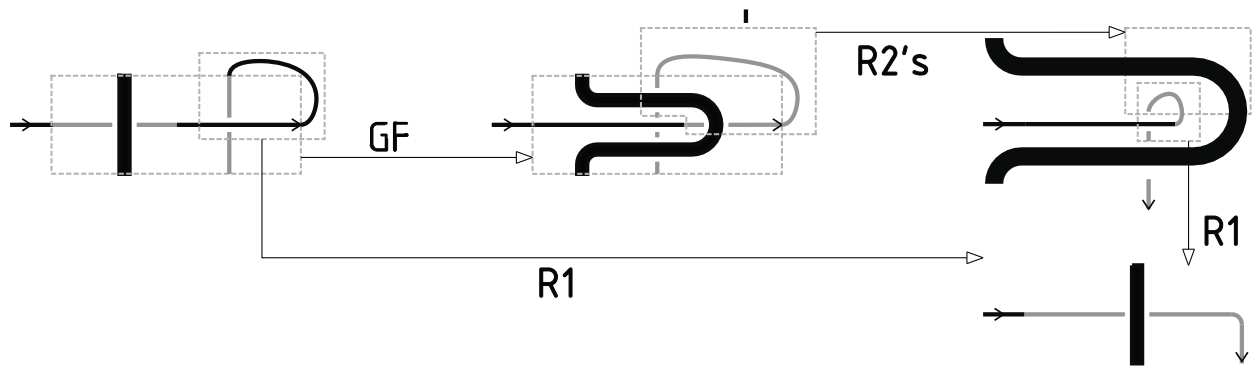}
\caption{Syzygy between $R1$- and $GF$- sorting moves.}
\label{fig:R1GFSyzygy}
\end{center}\end{figure}


\pagebreak
\subsection{Examples}
\label{sec:Examples}
Here are some examples of the sorting map applied to descending virtual long knots and pure tangles.
\begin{enumerate}%

\item The sorting map is applied to a generic framed descending virtual long knot diagram below:   where ``deform'' mean redraw the same virtual knot diagram on the plane,
\begin{center}
\includegraphics[width=6.5in]{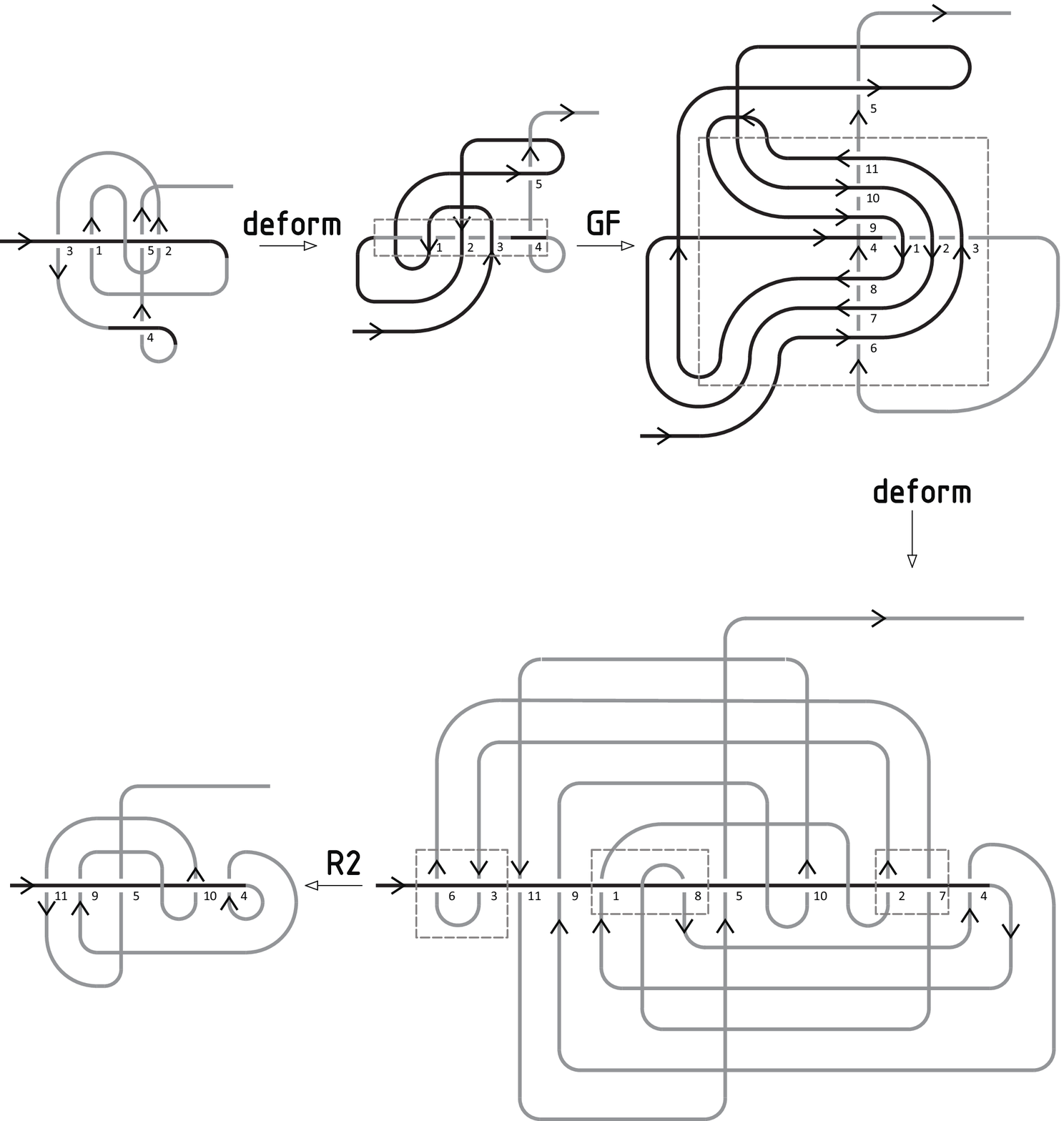}
\label{fig:SortExample}
\end{center}
Here the knot diagram is first ``deformed'' (or reimmersed on the plane) to show the one illegal interval (in the box in the second diagram), and then $GF$-sorted to remove the illegal interval, then deformed again to show the forbidden bigons, and finally $R2$-sorted to remove all bigons.  If the diagram represented an unframed virtual knot, then in this case an $R1$-sort can be used to arrive at the canonical form directly.  Notice that the canonical form obtained this way is the same as the one obtained by performing an $R1$-sort to the final diagram in the sequence above.

%

\item Two descending virtual pure tangle diagrams on three strands are sorted as follows:
\begin{center}
\includegraphics[width=6in]{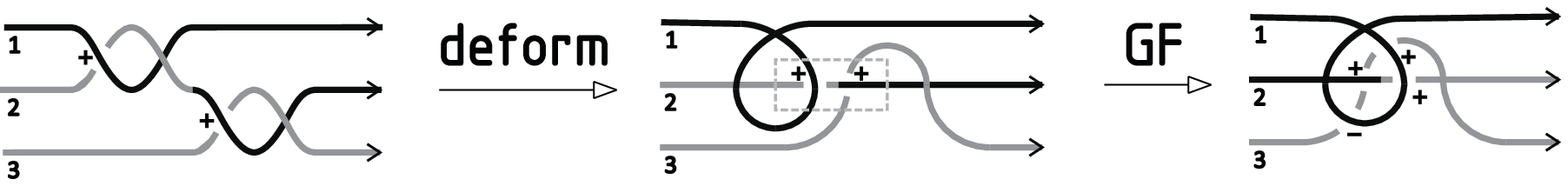}
\end{center}
\begin{center}
\includegraphics[width=6in]{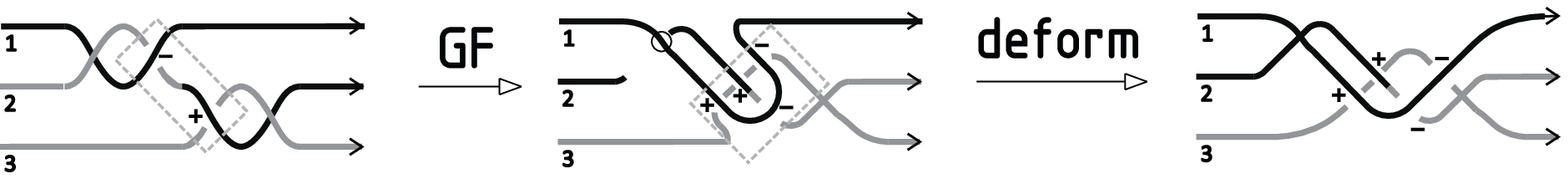}
\end{center}
Note that both starting diagrams are in ``braid-form,'' i.e. that as Gauss diagrams, the chords can be drawnparallel, but the canonical diagram for the first one does not remain in braid-form.

\begin{center}
\includegraphics[width=6in]{SortingTangleExample1.eps}
\end{center}
\begin{center}
\includegraphics[width=6in]{SortingTangleExample2.eps}
\end{center}
\end{enumerate}

\subsection{Remarks}
\begin{enumerate}
\item  For any descending virtual long knot, either of the two fundamental groups, defined by the Wirtinger presentation with base point either above or below the blackboard and which ignores the virtual crossings, is cyclic.  
 
 \item The classification of flat virtual pure tangles can be used as an invariant on virtual pure tangles as well as on virtual pure braids, presented by
$$v\mathcal{B}_n\,:=\,\langle \,\sigma_{ij} \,|\, \sigma_{ij} \sigma_{ik} \sigma_{jk} \,= \,\sigma_{jk} \sigma_{ik} \sigma_{ij}\, ,\, \sigma_{ij} \sigma_{kl} \,= \,\sigma_{kl} \sigma_{ij},\,\,\, 1\,<\,i,j,k,l\,\leq\,n\,\,\rangle$$
where $\sigma_{ij}$ can be represented by the positive crossing with strand $i$ over strand $j$ as follows:
\begin{center}
\includegraphics[height=0.4in]{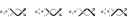}.
\end{center}
The virtual braid group on $n$ strands has an obvious map into the virtual pure tangles on $n$ strands.
\end{enumerate}

\begin{conjecture}
We conjecture that the flat virtual pure braid group on $n$ strand, the quotient of the $v\mathcal{B}_n$ by the flatness relation $\sigma_{ij}=\sigma_{ji}^{-1}$, injects into flat virtual pure tangles on $n$ strand.  If this is true, the classification above gives normal forms for the group which are not in terms of the alphabets in the presentation of the group.
\end{conjecture}

%


\end{document}